\documentclass[twoside,11pt,leqno]{article}
\usepackage{graphicx}
\def\reals{\hbox{\sl I\kern-.18em R \kern-.3em}}
\def\ints{\hbox{\sl Z\kern-.4em Z \kern-.3em}}
\def\nats{\hbox{\sl I\kern-.16em N \kern.05em}}
\def\rats{\hbox{\sl Q \kern-.83em\vrule height.59em depth0em \kern.87em}}
\def\complexes{\hbox{\sl\kern.50em I\kern-.50em C \kern.05em}}
\newtheorem{theorem}{Theorem}

\def\pf {{\bf Proof:} \ }
\def\endpf{$\|$ \medskip}

\def\be{\begin{enumerate}}
\def\ee{\end{enumerate}}
\def\bi{\begin{itemize}}
\def\ei{\end{itemize}}
\def\bs{\bigskip}
\def\ms{\medskip}

\def\cA{{\cal A}}
\def\cB{{\cal B}}

\def\cK{{\cal K}}
\def\cL{{\cal L}}

\def\cT{{\cal T}}

\def\cX{{\cal X}}

\def\ta{{\kern.35em{\cal C}\kern-.35em{\cal A}}}
\def\pa{{\kern.35em{\cal P}\kern-.35em{\cal A}}}

\def\bs{\bigskip}
\def\ms{\medskip}

\begin{document}

\markboth{Joan S. Birman}{ Braids, Knots and Contact structures}
\title{Braids, Knots and Contact Structures}
\author{ Joan S. Birman\\
Barnard College, Columbia University,\\
 Department of Mathematics, \\
 2990 Broadway, \\
 New York, New York 10027.\\
 jb\@math.columbia.edu}
\date{ February 8, 2004}
\maketitle

\begin{abstract}
These notes supplement my planned talk on Feb 19, 2004, at the First East Asian School of Knots and Related Topics, Seoul, South Korea.  I will review aspects of the interconnections between braids, knots and contact structures on $\reals^3$. I will discuss  my recent work with William Menasco \cite{BM-stab-I} and \cite{BM-stab-II}, where we prove that there are distinct transversal knot types in $\reals^3$ having the same topological knot type and the same Bennequin invariant.
\end{abstract}

 A \underline{knot} in oriented $\reals^3$
is the image $X$ of an oriented circle $S^1$ under a smooth embedding $e:S^1 \to \reals^3$.
Viewing $S^3$ as $\reals^3 \cup {\infty}$, we also can think of $X$ as being a knot in $S^3$.
The  \underline{knot type} $\cX$ of $X$ is its equivalence
class under smooth isotopy of the pair $(X,S^3)$.  \bs
 
  Let ${\bf A}$ be the $z$ axis in $\reals^3$, with standard cylindrical coordinates $(\rho,\theta, z)$ and let ${\bf H}$ be the collection of all half-planes $H_\theta$ through ${\bf A}$. The pair $({\bf A},{\bf H})$ defines the  \underline{standard} \underline{braid} \underline{structure} on $\reals^3$.  See the left sketch in Figure 1.
 \begin{figure}[htpb!]
 \label{figure:braid-contact}
\centerline {\includegraphics[scale=.7, bb=53 495 578 698]{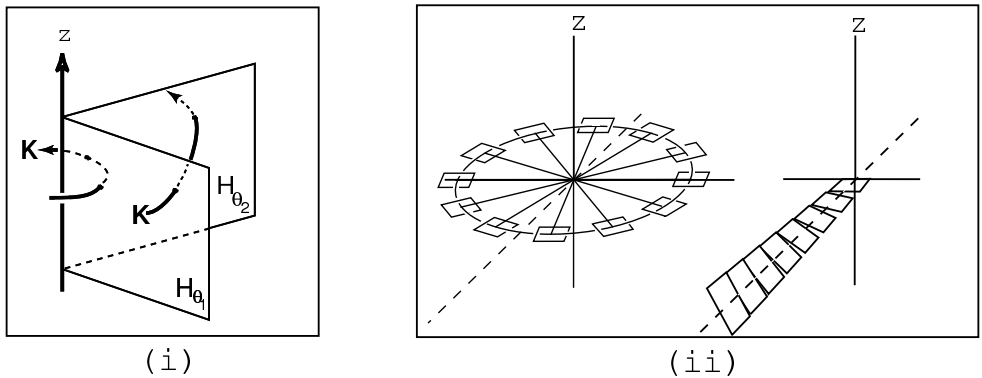}}
\caption{ (i) The standard braid structure in $\reals^3$; (ii)The standard polar contact structure on $\reals^3$}
\end{figure}
  Using the same cylindrical coordinates, let $\alpha$ be the 1-form $\alpha = \rho^2 d\theta + dz$.  The kernel $\xi$  of $\alpha$ defines  a contact structure on $\reals^3$.  The associated plane field is sketched on the right in Figure 1.  The family of 2-planes that define $\xi$ twist anticlockwise as one moves along the $x$ axis  from $0$ to $\infty$.  The family is invariant under rotation of 3-space about the $z$ axis and under translation of 3-space along rays parallel to the $z$ axis.  
The plane field has the property that $\alpha\wedge d\alpha$ is nowhere zero, which means that it is totally non-integrable, i.e. there is no surface in $\reals^3$ which is everywhere tangent to the 2-planes of $\xi$ in any neighborhood of any point in $\reals^3$. Intuitively, the twisting prevents tangencies between  surfaces and the collection of 2-planes. It is generic in the sense that, if $p$ is a point in any contact 3-manifold $M^3$, then in every neighborhood of $p$ in $M^3$ the contact structure is locally like the one we depicted in Figure 1.   Globally, in $\reals^3$ the total twist angle must be an odd multiple of $\pi$ for the contact structure to extend to $S^3$, and we distinguish the two cases by the symbols $\xi_\pi$ and $\xi_{>\pi}$.   The former (for reasons that will become clear shortly) is known as the \underline{standard} (polar) \underline{contact structure} and all of the latter are `overtwisted' contact structures. 
\ms

Let $K$ be a knot (for simplicity we restrict to knots here, but everything works equally well for links)  which is parametrized by cylindrical coordinates $(\rho(t), \theta(t), z(t))$, where $t \in [0,2\pi]$.    Then $K$ is a closed braid if $\rho(t)>0$ and $d\theta/dt > 0$ for all  $t$.  On the other hand, $K$ is a \underline{Legendrian} (resp. \underline{transversal}) knot if it is everywhere (resp. nowhere) tangent to the 2-planes of $\xi$.  In the Legendrian case this means that on $K$ we have  $d\theta/dt = (-1/\rho^2)(dz/dt)$.  In the transversal case we require that $d\theta/dt > (-1/\rho^2)(dz/dt)$ at every point of $K(t)$.  It was proved by Alexander in 1925 that every knot could be changed to a closed braid. Sixty years after Alexander proved his theorem, Bennequin adapted Alexander's proof to the setting of transversal knots in \cite{Bennequin}:, where he showed that  every transversal knot is isotopic, through transversal knots, to a closed braid.  \ms

\begin{figure}[htpb]
\centerline{\includegraphics[scale=.7]{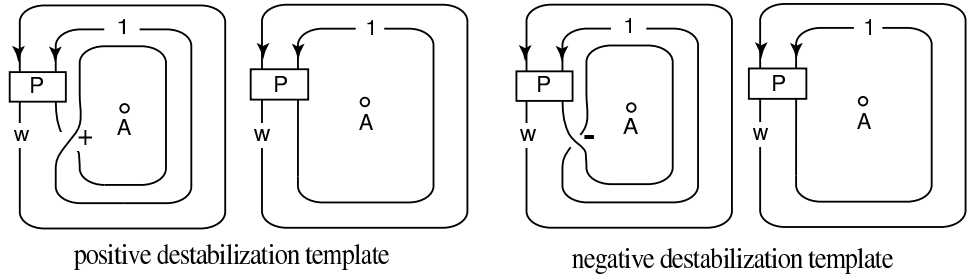}}
\caption{The two destabilization templates}
\label{figure:stab-destab}
\end{figure}

Closed braid representations of $\cX$ are not unique, and Markov's well-known 
theorem \cite{Birman1974} asserts that any two are related by a finite
sequence of elementary moves. 
One of the moves is \underline{braid isotopy}
\index{braid isotopy}, 
by which we mean an isotopy of the pair ($X, \reals^3\setminus{\bf A}$) which preserves the condition
that $X$ is transverse to the fibers of ${\bf H}$.
The  other two moves are mutually inverse, and are illustrated in Figure
\ref{figure:stab-destab}. 
Both take closed braids to closed braids. We call them
\underline{destabilization}  and
\underline{stabilization}, where the former decreases braid
index by one and the latter increases it by one. The weights, e.g.  $w$, that are attached to some of the strands  denote that many `parallel' strands, where parallel means in the framing defined by the given projection. The braid inside the box which is labeled
$P$ is an arbitrary $(w+1)$-braid. Later, it will be necessary to distinguish between
positive and negative destabilizations, so we illustrate both now.  The term
`templates', mentioned in the caption for Figure \ref{figure:stab-destab}, will be explained shortly.

\begin{theorem}
\label{theorem:MT}
{\bf (Markov (MT) \cite{Markov1935}):} Let $X_+,X_-$ be closed
braid representatives of the same oriented link type $\cX$ in oriented
3-space.  Then there exists a sequence of closed braid representatives
of $\cX$:
$$ X_+ = X_1 \to X_2 \to \cdots \to X_r = X_- $$
such that, up to braid isotopy, each $X_{i+1}$ is obtained from $X_i$ by a single stabilization or
destabilization. 
\end{theorem}

It is easy to find examples of subsequences $X_j \to \cdots \to X_{j+k}$ of (1) in Theorem
\ref{theorem:MT} such that
$b(X_j) = b(X_{j+k})$, but  $X_j$ and $X_{j+k}$
 are not braid isotopic.\ms

Seventy years after Markov's theorem was announced (a proof was not published until many years later, although at least 5 essentially different proofs exist today),  Orevkov and Shevchishin proved a version of the MT which holds in the transversal setting: 

\begin{theorem}
\label{theorem:TMT}
{\bf (Orevkov and Shevchishin  \cite{O-S}):} Let $TX_+,TX_-$ be closed
braid representatives of the same oriented link type $\cX$ in oriented
3-space.  Then there exists a sequence of closed braid representatives
of $\cT\cX$:
$$ TX_+ = TX_1 \to TX_2 \to \cdots \to TX_r = TX_1 $$
such that, up to braid isotopy, each $TX_{i+1}$ is obtained from $TX_i$ by a single positive stabilization or destabilization. 
\end{theorem}
Is the Transverse Markov Theorem  really different from the Markov Theorem?  Are there transversal knots which are isotopic as topological knots but are not transversally isotopic?  To answer this question we take a small detour and review the contributions of Bennequin in \cite{Bennequin}. \ms

Why did topologists begin to think about contact structures, and analysts begin to think about knots?  While we might wish that analysts suddenly became overwhelmed with the beauty of knots, there was a more specific and focussed reason.   At the time that Bennequin did his foundational work in \cite{Bennequin} it was not known  whether a 3-manifold could support more than one isotopy class of contact structures.  Bennequin answered this question in the affirmative, in the case of contact structures on $\reals^3$ or $S^3$ which were known to be homotopic to the standard one.   His  tool for answering it was highly original, and it had to do with braids and knots.  Let $TK$ be a transversal knot. Let $\cT\cK$ be its transversal knot type, i.e. its knot type under transversal isotopies, and let $[\cT\cK]_{top}$ be its topological knot type.  Choose a representative $TK$ of $\cT\cK$, which (by Bennequin's transversal version of Alexander's theorem) is always possible.  Choose a Seifert surface $F$ of minimal genus, with $TK = \partial F$.  Bennequin studied the foliation of $F$ which is induced by the intersections of $F$ with the plane field determined by $\xi_\pi$.     Let $n(TK)$ be the braid index and let $e(TK)$ be the algebraic crossing number of a projection of $TK$ onto the plane $z=0$.  Both can be determined from the foliation.  Bennequin found an invariant of $\cT\cK$, given by the formula $\beta(\cT\cK) = e(TK) - n(TK)$.  Of course if he had known Theorem \ref{theorem:TMT}  the proof that $\beta(\cT\cK)$ is an invariant of $\cT\cX$ would have been trivial, but he did not have that tool.  He then showed a little bit more: he showed that $\beta(\cT\cK)$ is bounded above by the negative of the Euler characteristic of $F$ in $\xi_\pi$.  He then showed that this bound fails in one of the contact structures $\xi_{>\pi}$. In this way he proved that the contact structures  $\xi_{>\pi}$ cannot be isotopic to $\xi_\pi$. To knot theorists, his proof should seem intuitively natural, because the invariant $\beta(\cT\cK)$ is a self-linking number of a representative $TX\in \cT\cX$ (the sense of push-off being determined by $\xi$), and the more twisting there is the higher this number can be.  For an explanation of the self-linking, and lots more about Legendrian and transversal knots we refer the reader to John Etnyre's excellent review article \cite{Etnyre2003}. The basic idea is that $TX$ bounds a Seifert surface, and this Seifert surface is foliated by  the plane field associated to $\xi$. Call this foliation the \underline{characteristic foliation}. Near the boundary, the characteristic foliation is transverse to the boundary. The Bennequin invariant is the linking number of $TX, TX'$, where $TX'$ is is a copy of $TX$, obtained by pushing $TX$ off itself onto $F$, using the direction determine by the characteristic foliation of $F$. \ms

Bennequin's paper was truly important.   Shortly after it was written Eliashberg showed in \cite{Eliashberg} that the phenomenon of an infinite sequence of contact structures related to a single one of minimal twist angle occured generically in every 3-manifold, and introduced the term `tight' and `overtwisted' to distinguish the two cases.   Here too, there is a reason that will seem natural to topologists. In 2003 Giroux proved \cite{Giroux} that every contact structure on every closed, orientable 3-manifold $M^3$ can be obtained in the following way: Represent $M^3$ as a branched covering space of $S^3$, branched over a knot or link, and lift the standard and overtwisted contact structures on $S^3$  to $M^3$.  \ms  
 
Returning to knot theory, the invariant $\beta(\cT\cX)$ allows us to answer a fundamental question: is the equivalence relation on knots that is defined by transversal isotopy really different from the equivalence relation defined by topological isotopy?
\begin{theorem}
\label{theorem:transversal isotopy}
{\bf (Bennequin \cite{Bennequin}):}  There are infinitely many distinct transversal knot types associated to each topological knot type.
\end{theorem}

\pf  Choose a transversal knot type $\cT\cK$ and a closed braid representative $TX_0$.  Stabilizing the closed braid $TX_0$ once negatively, we obtain the transverse closed braid $TX_1$, with $e(TX_1) = e(TX_0)-1$ and $n(TX_1) = n(TX_0) + 1$, so that $\beta(TX_1) = \beta(TX_0) - 2$. Iterating, we obtain transverse closed braids $TX_2, TX_3, \dots$, defining transverse knot types $\cT\cX_1, \cT\cX_2, \cT\cX_3, \cdots$,  and no two have the same Bennequin invariant. Since stabilization does not change the topological knot type, the assertion follows.  \endpf

This brings us to the research that is the main goal of this review article.  It is an outstanding open problem to find computable invariants of $\cT\cX$ which are not determined by $[\cT\cX]_{top}$ and $\beta(\cT\cX)$.   A hint that the problem
might turn out to be quite subtle was in the paper \cite{F-T} by Fuchs and Tabachnikov, who proved that while ragbags
filled with polynomial and finite type invariants of transversal knot types $\cT\cX$ exist, based upon the work of
V.I. Arnold in \cite{Arnold}, they are all determined by
$[\cT\cX]_{top}$ and $\beta(\cT\cX)$. Thus, the seemingly new invariants that many people had discovered by using
Arnold's ideas were just a fancy way of encoding $[\cT\cX]_{top}$ and $\beta(\cT\cX)$. \ms

This leads naturally to a question:  Are there computable invariants of transversal knots which are {\it not} determined by  $[\cT\cX]_{top}$ and $\beta(\cT\cX)$?  A similar question arises in the setting of Legendrian knots.  Each Legendrian knot $\cL\cX$ determines a topological knot type $[\cL\cX]_{top}$, and just as in the transverse case it is an invariant of the Legendrian knot type. There are also two numerical invariants of $\cL\cX$: the Thurston-Bennequin invariant $tb(\cL\cX)$ (a self-linking number) and the Maslov index $M(\cL\cX)$ (a rotation number).   So until a few years ago the same question existed in the Legendrian setting. But the Legendrian case has been settled.  

\begin{theorem} {\bf (Chekanov  \cite{Chekanov}):} There exist distinct Legendian knot types  which have the same topological knot type$[\cL\cX]_{top}$, and also the same Thurston-Bennequin invariant $tb(\cL\cX)$ and Maslov index $M(\cL\cX)$. 
\end{theorem}

\begin{figure}[htpb]
\centerline{\includegraphics[scale=.7, bb = 124 573 460 712]{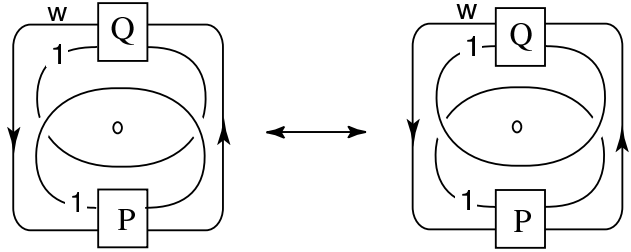}}
\caption{The exchange move template}
\label{figure:exchange}
\end{figure}

The analogous result for transversal knots proved to be quite difficult, so to begin to understand whether something could be done via braid theory Birman and Wrinkle asked an easier question: are there {\it are} transversal knot types which {\it are determined}  by their topological knot type and Bennequin number?  This question lead to a definition
in \cite{Birman-Wrinkle}: a transversal knot type $\cT\cX$ is \underline{transversally simple} if it is determined by
$[\cT\cX]_{top}$ and $\beta(\cT\cX)$. We asked: are there transversally simple knots?  The manuscript \cite{Birman-Wrinkle} gives a purely topological (in fact braid-theoretic) criterion which enables one to answer the question affirmatively, adding one more piece of evidence that topology and analysis walk hand in hand.  To explain what we did, note that there is no loss in generality in working in the setting of closed braids.  It will be convenient to introduce a new move that takes closed braid to closed braids, the \underline{exchange move}.  See Figure \ref{figure:exchange}.   It is easy to prove that the exchange move can be realized as a transversal isotopy between transversal closed braids, so that 
while (by Theorem \ref{theorem:TMT}) it must be a consequence of braid isotopy and positive stabilizations and destabilizations, there is no harm in adding the exchange move to the bag of tools that one has in simplifying transversal closed braid representatives of transversal knots.
A topological knot or link type  $\cX$ is said to be \underline{exchange reducible} if an arbitrary closed braid representative $X$ of $\cX$ can be changed to an arbitrary representative of minimum braid index  
by braid isotopy, positive and negative destabilizations and exchange moves.  We have:

\begin{theorem}  {\bf (Birman and Wrinkle \cite{Birman-Wrinkle}):}
\label{theorem:exchange-reducibility}
If a knot type $\cX$ is exchange-reducible, then any transversal knot type $\cT\cX$ which has $[\cT\cX]_{top} = \cX$  is transversally simple. 
\end{theorem}
This theorem was used to give a new proof of a theorem of Eliashberg \cite{Eliashberg1993}, which asserts that the unlink is transversally simple, and also (with the help of \cite{Menasco2001}) to prove the then-new result that most iterated torus knots are transversally simple.\ms

The rest of this review will be directed at explain the main result of \cite{BM-stab-II}, joint work of the author and W. Menasco:  
\begin{theorem}{\bf (Birman and Menasco \cite{BM-stab-II}):}
\label{theorem:negative flype examples}
There exist transversal knot types which are not transversally simple. Explicitly, the  transverse closed $3$-braids 
$TX_+ =  \sigma_1^5 \sigma_2^4 \sigma_1^6\sigma_2^{-1}$ and
  $TX_- = \sigma_1^5 \sigma_2^{-1} \sigma_1^6\sigma_2^4$
determine transverse knot types $\cT\cX_+, \cT\cX_-$ with $(\cT\cX_+)_{top} = (\cT\cX_-)_{top}$ and $\beta(\cT\cX_+) = \beta(\cT\cX_-)$, but  $\cT\cX_+ \not= \cT\cX_-$.
\end{theorem} 
We remark that the proof of Theorem \ref{theorem:negative flype examples} does not use a computable invariant, rather it is very indirect.  The problem of finding new computable invariants of transversal knot types remains open at this writing. 

Before we can describe the proof of Theorem \ref{theorem:negative flype examples}, we need to explain the Markov Theorem Without Stabilization (MTWS), established by the author and Menasco in \cite{BM-stab-I}.  As in the case of Theorems \ref{theorem:MT} and \ref{theorem:exchange-reducibility},   the moves that are needed for the MTWS will be  described in terms of pairs of `block-strand
diagrams' which we call `templates'.  Examples of block-strand diagrams were given in Figures \ref{figure:stab-destab} and \ref{figure:exchange}. Their important feature is that after an assignment of a braided tangle to each block, the block strand diagram becomes a closed braid.  Our block strand diagrams come in pairs. A pair of block-strand diagrams are called a \underline{template}, and the templates define the `moves' of the theorem. Here the important feature is that  the two block-strand diagrams in a template represent the same knot or link,  for \ {\it every} \  braiding assignment to the blocks. 

\begin{figure}[htpb!]

\centerline{\includegraphics[scale=.7, bb = 30 483 566 731]{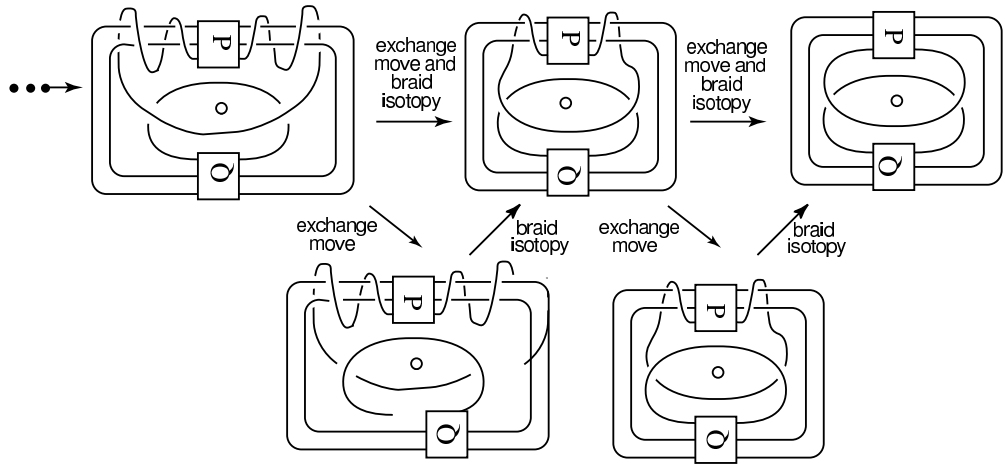}}

\caption{The exchange move can lead to infinitely many distinct conjugacy classes of $n$-braids representing the same oriented link type} \label{figure:exchiso}

\end{figure}

The  exchange move looks very harmless, but unfortunately it leads to Markov towers which produce infinitely many closed braid representatives of a knot or link, which (for almost all braiding assignments to the blocks $P$ and $Q$) can be shown to be in distinct braid isotopy (or conjugacy) classes .   See Figure \ref{figure:exchiso}.  
The phenomenon which is exhibited in Figure \ref{figure:exchiso} complicates the statement of the Markov Theorem Without Stabilization (MTWS), which is our next result:

\begin{theorem}{\bf (Birman and Menasco \cite{BM-stab-I}):}
\label{theorem:MTWS}
Let $\cB$ be the collection of all braid isotopy classes of closed braid representatives of oriented
knot and link types in oriented 3-space. Among these, consider the subcollection
$\cB(\cX)$ of representatives of a fixed link type $\cX$.  Among these, let $\cB_{min}(\cX)$ be the
subcollection of representatives whose braid index is equal to the braid index of $\cX$.  Choose any
$X_+\in \cB(\cX)$  and any $X_-\in \cB_{min}(\cX)$.  Then there is a
complexity function  which is associated
to $X_+,X_-$, and for each braid index $m$ a finite set $\cT(m)$ of templates is introduced, each
template determining a move which is non-increasing on braid index, such that the following hold:
First, there is are initial sequences which modify $X_- \to X_-'$ and $X_+ \to X_+'$:
$$X_- = X_-^1 \to \cdots\to X_-^p = X_-',   \ \ \ \ \ X_+ = X_+^1\to\dots \to X_+^q = X_+'  $$
Each passage $X_-^j \to X_- ^{j+1}$  is strictly complexity
reducing  and is realized by an exchange move, so that $b(X_-^j) = b(X_-^{j+1}) $.  These moves `unwind' $X_-$, if it is wound up as in the top right sketch in Figure \ref{figure:exchiso}.  Each passage
$X_+^j \to X_+^{j+1}$   is strictly  complexity-reducing and is realized by either an exchange move or a destabilization, so that $b(X_+^j) \geq b(X_+^{j+1})$.   Replacing $X_+$ with $X_+'$ and $X_-$ with $X_-'$, there is an additional
sequence which modifies $X_+'$, keeping $X_-'$ fixed:
$$ X_+'=X^q\to\cdots \to X^r = X_-' $$
Each passage $X^j\to X^{j+1}$ in this sequence  is also
strictly complexity-reducing.  It is realized by an exchange move,  destabilization, or one of the moves defined by a template
$\cT$ in the finite set $\cT(m)$, where $m=b(X_+)$. The inequality $b(X^j) \geq b(X^{j+1})$ holds
for each $j = q,\dots,r-1$ and so also for each $j=1,\dots,r-1$.
\end{theorem}

\begin{figure}[htpb!]
\label{figure:flype}

\centerline{\includegraphics[scale=.7, bb = 53 567 571 700]{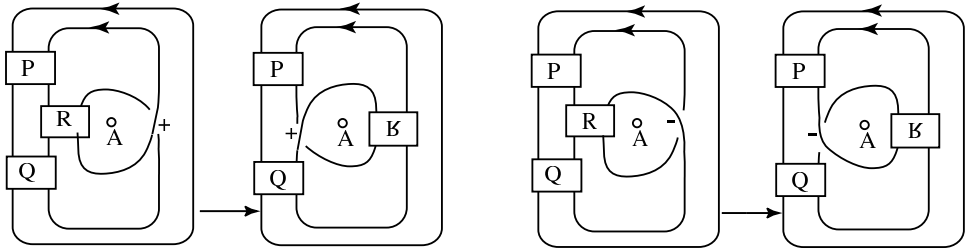}}

\caption{The two flype templates differ in the sense of the half-twist that realizes the isotopy. We call the left (resp. right) one positive (resp. negative).} 

\end{figure}

Figure \ref{figure:flype} shows two more examples of the templates of the MTWS, namely the two flype templates.  Many examples of more complicated templates are given in the manuscript \cite{BM-stab-I}.  In these more general templates the Markov towers are quite complicated, and so the isotopy that takes the left diagram  to the right one is often not obvious. \ms

The proof of Theorem \ref{theorem:MTWS} may be described {\it very} briefly as follows.    We are given two closed braids, $X_+$ and $X_-$, and an isotopy that takes $X_+$ to $X_-$.  The trace of the isotopy sweeps out an annulus, but in general it is not embedded.  The proof begins by showing that the given isotopy can be split into two parts, over which we have some control. An intermediate link $X_0$ which represents the same link type $\cX$ as $X_+$ and $X_-$  is constructed, such that the trace of the isotopy from $X_+$ to $X_0$ is an embedded annulus $\cA_+$.  Also the trace of the isotopy from $X_0$ to $X_-$ is a second embedded annulus $\cA_-$. The union of these two embedded annuli $\cT\cA = \cA_+\cup\cA_-$  is an immersed annulus, but its self-intersection set is controlled, and is a finite number of clasp arcs.   The main tool in the proof of Theorem \ref{theorem:MTWS} is the study of the foliation of the  immersed annulus $\cT\cA$, which is induced by its intersections with the half-planes of the braid structure that was illustrated earlier in the left sketch of Figure 1   We remarked earlier that as one moves away from the braid axis the braid foliations of a surface bounded by a knot or link resemble, in key ways, the foliation induced by the standard contact structure.  Therefore it should seem natural to the reader that  the MTWS plays a key role in the proof of Theorem \ref{theorem:negative flype examples}, which we discuss next.  This is yet another instance of the main theme of this little review: the close connections between the mathematics of braids and the mathematics of contact structures.\ms

\noindent {\bf Sketch of the proof of Theorem \ref{theorem:negative flype examples} \cite{BM-stab-II}}
First we show that the examples satisfy the conditions of the theorem.  The topological knot types defined by the closed 3-braids $ \sigma_1^5 \sigma_2^4 \sigma_1^6\sigma_2^{-1}$ and $\sigma_1^5 \sigma_2^{-1} \sigma_1^6\sigma_2^4$  coincide because they are carried by the block strand diagrams for the negative flype template
of Figure \ref{figure:flype}.  The  Bennequin invariant can be computed as the exponent sum of the braid word (14 in both cases) minus the braid index (3 in both cases).  So our examples have the required properties.  \ms

The hard part is the establishment of a special version of Theorem \ref{theorem:MTWS} which is applicable to the situation that we face.  Its special features are as follows: \ms

\noindent (1) Both $X_+$ and $X_-$ have braid index 3.   Since it is well known that exchange moves can be replaced by braid isotopy for 3-braids,   the first two sequences in Theorem \ref{theorem:MTWS} are vacuous, i.e. $X_\pm =X_\pm^\prime$.\ms

\noindent (2) With the restrictions in (1) above,  it is proved in \cite{BM-stab-II} that if $X_-$ and $X_+$ are transversal closed braids $TX_+$ and $TX_-$, then the isotopy that takes $TX_+$ to $TX_-$ may be assumed to be transversal.\ms

\noindent (3) Because of the special assumption, the templates that are needed, in the topological setting, can be enumerated explicitly: they are the positive and negative destabilization and flype templates. No others are needed.  \cite{K-L1999}  \ms

Suppose that a transversal isotopy exists from the transverse closed braid $TX_+$ to the transverse closed braid $TX_-$.  Then (2) above  tells us that there is a $3$-braid template that carries the braids $ \sigma_1^5 \sigma_2^4 \sigma_1^6\sigma_2^{-1}$ and $\sigma_1^5 \sigma_2^{-1} \sigma_1^6\sigma_2^4$.  This is the first key fact that we use from Theorem \ref{theorem:MTWS}.  Instead of having to consider {\it all} possible transversal isotopies from $TX_+$ to $TX_-$, we only need to consider ones that relate the left and right block-strand diagrams in one of the four 3-braid templates.  By (3) above, the braids in question are carried by either one of the two destabilization templates or one of the two flype templates.   If it was one of the destabilization templates, then the knots in question could be represented by 2 or 1-braids, i.e. they would be type $(2,n)$ torus knots or the unknot, however  an easy argument shows that the knots in Theorem \ref{theorem:negative flype examples} are neither type $(2,n)$ torus knots or the unknot.  The positive flype templates are ruled out in different way: We know that, topologically, our closed braids admit a negative flype, so if they are also carried by the positive flype template they admit flypes of both signs. But the manuscript \cite{K-L1999} gives conditions under which a closed 3-braid admits flypes of both signs, and the examples were chosen explicitly to rule out that possibility.  \ms

We are reduced to isotopies that are supported by the negative flype template.  We know that the obvious isotopy is not transversal, but maybe there is some {\it other} isotopy which is transversal. Here we use a key fact about the definition of a template (and this is a second very strong
aspect of the MTWS).  If such a transversal  isotopy exists, then it exists for every knot or  link defined by a fixed choice of braiding assignments to the blocks.  Choose the braiding assignments
$\sigma_1^3, \sigma_2^4, \sigma_1^{-5}$ to the blocks
$P,R,Q$.  This braiding assignment gives a 2-component link
$L_1\sqcup L_2$ which has two distinct isotopy classes of closed 3-braid representatives. If $L_1$ is the
component associated to the left strand entering the block $P$, then
$\beta(L_1)=-1$ and $\beta(L_2) = -3$ before the flype, but after the flype the
representative will be $\sigma_1^3 \sigma_2^{-1} \sigma_1^{-5} \sigma_2^{4}$, with 
$\beta(L_1)=-3$ and $\beta(L_2) = -1$. However, by Eliashberg's isotopy extension theorem (Proposition
2.1.2 of \cite{Eliashberg1993}) a  transversal isotopy of a knot/link extends to an ambient
transversal isotopy  of the 3-sphere.  Any transversal isotopy of $L_1 \sqcup L_2$
must preserve the $\beta$-invariants of the components. It follows that no such transversal isotopy exists, a contradiction of our assumption that $TX_+$ and $TX_-$ are transversally isotopic. \endpf

\noindent {\bf Remark:} Other examples of a similar nature were discovered by Etnyre and Honda \cite{E-H2003} after the proof of Theorem \ref{theorem:negative flype examples} was posted on the arXiv.  Their methods are very different from the proof that we just described (being based on contact theory techniques rather than topological techniques), but are equally indirect. They do not produce explicit examples, rather they present a bag of pairs of transverse knots and prove that at least one pair in the bag exists with the properties given by Theorem \ref{theorem:negative flype examples}. \ms

 \end{document}